\documentclass[fleqn]{article}

\newcommand{\titel} {Subvarieties  of  the variety of meadows}

\usepackage{graphicx,version,stmaryrd,amssymb,amsmath,amsthm,url}
\usepackage{graphics}
\usepackage[usenames,dvipsnames]{color}

\newtheorem{theorem}{Theorem}
  
\newtheorem{proposition}{Proposition}  
  
\newtheorem{corollary}{Corollary}

\newtheorem{definition}{Definition}  
\theoremstyle{definition}

\newcommand{\EFR}{\mathsf{EFR}}
\newcommand{\Inv}{\mathsf{Inv}}
\newcommand{\IL}{\mathsf{IL}}
\newcommand{\sg}{\operatorname{{\mathbf s}}}
\newcommand{\Md}{\mathsf{Md}}
\newcommand{\SA}{\mathsf{Signs}}
\newcommand{\SR}{\mathsf{Square roots}}

\title{\titel}

\author{
	Jan A.\ Bergstra \& Inge Bethke \\
	\\
	 {\small
	  Informatics Institute,
	  University of Amsterdam}\\
	{\small Email: \url{j.a.bergstra@uva.nl}, \url{i.bethke@uva.nl}
	}
}
\date{}

\begin{document}
\maketitle

\begin{abstract}
Meadows|commutative rings equipped with a total inversion operation|can be axiomatized by purely equational means.
We study subvarieties of the variety of meadows obtained by extending the equational theory and expanding the signature.
\end{abstract}

\section{Introduction}

\begin{table}[t]
\centering
\hrule
\begin{align*}
	(x+y)+z &= x + (y + z)& (1)\\
	x+y     &= y+x& (2)\\
	x+0     &= x& (3)\\
	x+(-x)  &= 0& (4)\\
	(x \cdot y) \cdot  z &= x \cdot  (y \cdot  z)& (5)\\
	x \cdot  y &= y \cdot  x& (6)\\
	1\cdot x &= x & (7)\\
	x \cdot  (y + z) &= x \cdot  y + x \cdot  z& (8)\\
	(x^{-1})^{-1} &= x& (9) \\
	x \cdot (x \cdot x^{-1}) &= x& (10)
\end{align*}
\hrule
\caption{The set $\Md$ of axioms for meadows}
\label{meadow}
\end{table}
\begin{table}[t]
\centering
\hrule
\begin{align*}
	\forall x \exists y \ x&=x\cdot y\cdot x& (11)\\
\end{align*}
\hrule
\caption{The axiom for the weak inverse}
\label{win}
\end{table}

In \cite{vN36}, von Neumann introduced his regular rings during the study of von Neumann algebras and continuous geometry. 
A \emph{von Neumann regular ring} is an algebra $\langle R, +, \cdot, -, 0,1\rangle$ satisfying the axioms (1)--(5), (7), (8), (11) of Table \ref{meadow} and \ref{win}. In general, $y$ in (11) is not uniquely determined by $x$. However, if the underlying ring is commutative then the weak inverse is unique. 
\emph{Meadows} were introduced in \cite{BT2007} as commutative rings equipped with a total inversion operation. A meadow is an algebra $\langle R, +, \cdot, -, \ ^{-1}, 0,1\rangle$ satisfying the axioms (1)--(10) in Table \ref{meadow}. 
The paradigmatic  infinite meadows are the rational, real, and complex numbers each equipped with a total multiplicative inversion operation where $0^{-1}=0$. Recently, we found that  equational theories of  fields in which $0^{-1}=0$ holds were already introduced by Komori \cite{Kom75} in a report from 1975 under the notion of \emph{desirable pseudo-fields}.
Every commutative von Neumann regular ring can be obtained as a reduct of a meadow, and conversely, every meadow can be obtained from a commutative von Neumann regular ring by defining the inversion operation. Thus, although the signatures of meadows and commutative von Neumann regular rings differ, they form the same category of objects and natural morphisms.

From the perspective of universal algebra, the situation is quite different. If $K$  is a class of algebras of the same signature, then $K$ is a \emph{variety} if $K$ is closed under subalgebras, homomorphic images and direct products of nonempty families. One of the most celebrated theorems of Birkhoff says that the \emph{equational} classes of algebras|algebras axiomatized by identities|are precisely the varieties. Hence, the class of meadows  is a variety. In contrast, the class of von Neumann regular rings is merely an \emph{elementary} class|a class axiomatized by a set of first-order formulas|since it is not closed under subalgebras: e.g.\ the rational numbers $\mathbb{Q}$ form a von Neumann regular ring but its subalgebra $\mathbb{Z}$ of integers does not have weak inverses (except for 0 and 1).  In \cite{H89}, the concept of an \emph{existence variety} is introduced which applies to (commutative) von Neumann regular rings.  Existence  varieties require  only the closure of a class of algebras under elementary substructures, homomorphic images and direct products but also enjoy a Birkhoff-style theorem.

{\mbox{}A challenging question is whether there exist finite complete axiomatizations for the meadows of the rational, real, and complex numbers.
\cite{BBP15} gives an affirmative answer in the case of the expansion of the meadow of the real numbers with a sign function. 
In this paper we try to attack the problem for the meadow of the rational numbers by studying the structure of the subvarieties of meadows.} We presume a general familiarity with universal algebra and initial algebra semantics. Some  references to various aspects of this subject are \cite{BS81,GTWW77,Gr79,MMT87,W92}. In Section 2, we determine the simple meadows|a special kind of subdirectly irreducible algebras having no nontrivial congruence relations. It turns out that the simple meadows are exactly the building blocks of meadows. In Section 3, we study subvarieties of meadows obtained by  extending the theory.  In Section 4, we expand meadows with additional signature. In special cases we can give an upper bound  of the number of equations needed to axiomatize finitely based theories. In particular, we consider the expansion of meadows with a sign function. In Theorem \ref{sign} we prove that the sign function has an equational specification in the signature of the meadows.  In Section 5, we discuss presentations of meadows and specializations of meadow varieties. 
Section 6 ends the paper with open questions.

\section{The cancellation  meadows}
In the sequel, we shall write $\frak{Alg}(T)$ for the class of algebras axiomatized by the  theory $T$. Thus e.g.\ $\frak{Alg}(\Md)$ denotes the variety of meadows. Moreover, if $\mathcal{F}$ is a field, we shall write $\mathcal{F}_0$ for the expansion of $\mathcal{F}$ where  the inversion operation is completed  by $0^{-1}=0$.  Thus e.g.\ $\mathbb{Q}_0$ denotes the expansion field of the rational numbers where  the inversion operation is completed  by $0^{-1}=0$. 

The `meaning' of an equational class $V$ is often taken as the \emph{initial algebra} of that class, i.e.\  the algebra $\mathcal{I}_V\in V$ which has the special property that for every algebra $\mathcal{A}\in V$  there exists a unique homomorphism $h:\mathcal{I}_V \rightarrow \mathcal{A}$. The initial algebra always exists, is unique up to isomorphism, and can be constructed from the closed term algebra by  dividing out over provable equality. In \cite{BR10}, it is shown that the \emph{initial meadow}|the initial algebra in $\frak{Alg}(\Md)$|can also be represented as the minimal subalgebra of the direct product of the expansions of the finite prime fields, i.e.\  the finite meadows of the form $(\mathbb{Z}/p\mathbb{Z})_0$ with $p$ a prime number. Dually, an algebra $\mathcal{A}$ is \emph{final} in $V$ if for every $\mathcal{B}\in V$, there is a unique homomorphism $h: \mathcal{B} \rightarrow \mathcal{A}$. The \emph{final meadow}|the final algebra in $\frak{Alg}(\Md)$|is the trivial, one-element algebra.

A special class of meadows are the so-called \emph{cancellation} meadows|meadows that satisfy the \emph{Inverse Law}
\[
x\neq 0 \longrightarrow x \cdot x^{-1}=1. \ \ \ \ (\IL)
\]
Cancellation meadows and expansion fields $\mathcal{F}_0$ form the same class $\frak{Alg}(\Md + \IL)$|a class that is not closed under products: e.g.\  $(\mathbb{Z}/2\mathbb{Z})_0 \times (\mathbb{Z}/3\mathbb{Z})_0 \not \in \frak{Alg}(\Md + \IL)$ since $\langle 0,1\rangle \neq \langle 0,0\rangle $ but $\langle 0,1\rangle \cdot \langle 0,1\rangle ^{-1}=\langle 0,1\rangle \neq \langle 1,1\rangle$.
The class of cancellation meadows is therefore not a (sub)variety of meadows. In particular, cancellation meadows cannot be axiomatized by purely equational means.

An algebra $\mathcal{A}$ is a \emph{subdirect product} of an indexed family $(\mathcal{A}_i)_{i\in I}$ if $\mathcal{A}$ is a subalgebra of the direct  product $\prod_{i\in I} \mathcal{A}_i$ and $\pi_i(\mathcal{A})=\mathcal{A}_i$, for every $i\in I$. If $\mathcal{A}$ is an algebra, we denote by $Con(\mathcal{A})$ its lattice of congruences. Recall that $\mathcal{A}$ is \emph{subdirectly irreducible} if and only if $Con(\mathcal{A})-\{\Delta\}$|with $\Delta$ the diagonal|has a least element. The subdirect representation theorem of universal algebra states that every algebra is a subdirect product of subdirectly irreducible algebras. 
In \cite{BR10} it is proved that every subdirectly irreducible meadow is a cancellation meadow. We therefore have the following theorem.
\begin{theorem}\label{cancel}
Every meadow is a subdirect product of cancellation meadows. 
\end{theorem}
So cancellation meadows form the building blocks of meadows. It follows that the smallest variety of meadows containing all cancellation meadows is the entire class of meadows.
\begin{theorem}\label{cancellation}
Let $V$ be a variety with $ \frak{Alg}(\Md + \IL)\subseteq V\subseteq  \frak{Alg}(\Md)$. Then $V= \frak{Alg}(\Md)$.
\end{theorem}
{\bf Proof}: $V$ contains every meadow, since every meadow is a subdirect product of cancellation meadows and $V$ is closed under products and subalgebras. \hfill $\Box$
\begin{definition}
An algebra $\mathcal{A}$ is \emph{simple} if the only congruences on its carrier set $A$ are the diagonal $\Delta$ and the all relation $A \times A$.
\end{definition}
It is not hard to see that an algebra is simple if and only if it has no proper homomorphic images. Some textbooks require that a simple algebra be nontrivial. For our development we find the discussion smoother by admitting the  one-element meadow.
\begin{theorem}
$\frak{Alg}(\Md + \IL)$ is the class of simple meadows. 
\end{theorem}
{\bf Proof:} Let $\mathcal{M}$ be a cancellation meadow. Suppose the congruence  $\Theta$ on its carrier $M$ is neither the diagonal nor the all relation. Then we can pick $a_0, a_1, a_2, a_3\in M$ such that $a_0\neq a_1$, $\langle a_0,a_1\rangle \in \Theta$ and $\langle a_2,a_3\rangle\not\in \Theta$. Since $\Theta$ is a congruence we have $\langle a_0\cdot (a_2-a_3), a_1\cdot (a_2-a_3)\rangle\in \Theta$ and hence $\langle a_0\cdot a_2 - a_1 \cdot a_2, a_0\cdot a_3 -a_1\cdot a_3\rangle\in \Theta$. Therefore 
$\langle a_2\cdot (a_0 -a_1), a_3\cdot (a_0 -a_1)\rangle\in \Theta$ and thus $\langle a_2\cdot (a_0 -a_1)\cdot (a_0 -a_1)^{-1}, a_3(a_0 -a_1)\cdot (a_0 -a_1)^{-1}\rangle\in \Theta$. So $\langle a_2,a_3\rangle\in \Theta$ by ($\IL$).This yields a contradiction. Therefore $\Theta$ must be either the diagonal or the all relation, i.e.\  $\mathcal{M}$ is simple.

Conversely, let $\mathcal{M}$ be simple meadow. If $\mathcal{M}$ is final, then it is a cancellation meadow. Otherwise, $\mathcal{M}$ is a nontrivial subdirect product of cancellation meadows. Thus there is a projection $\pi$ with $\pi(\mathcal{M})$ a nontrivial cancellation meadow. Now, since $ker(\pi)$ is a congruence on the carrier $M$ of $\mathcal{M}$, $ker(\pi)=\{\langle a,a'\rangle\in A\times A\mid \pi(a)=\pi(a')\}$ must be the diagonal. Hence $\mathcal{M} \cong\pi(\mathcal{M})$. Hence $\mathcal{M}$ is a cancellation meadow.  \hfill $\Box$

\section{Subvarieties of meadows}
A variety is \emph{trivial} if it consists of the final one-element algebra only, otherwise it is \emph{nontrivial}. A sublclass  $W$ of a variety $V$ which is also a variety is called a \emph{subvariety}. Subvarieties of $\frak{Alg}(\Md)$ arise when equations are added to $\Md$.  In this section we study initial algebra specifications of $\mathbb{Q}_0$.

\begin{definition}
For $n\in \mathbb{N}$, we define the numeral $\underline{n}$ by $\underline{0} =0$ and $\underline{n+1}=\underline{n}+1$, and  for $P$|the set of prime numbers|we put $\Inv_P= \{\underline{p}\cdot \underline{p}^{-1}=1 \mid p \in P\}$.
\end{definition}
A necessary and sufficient condition for  an initial algebra specification of $\mathbb{Q}_0$ is given in the following theorem.
\begin{theorem}\label{initialspec}
Let $E$ be a set of meadow equations. If  $\frak{Alg}(\Md +E )$ is nontrivial, then
\[
\mathbb{Q}_0 \cong \mathcal{I}_{\frak{Alg}(\Md +E )}  \text{ if and only if  for all prime numbers $p$, } (\mathbb{Z}/p\mathbb{Z})_0 \not \models E.
\]
\end{theorem}
{\bf Proof}: Assume $\mathbb{Q}_0 \cong  \mathcal{I}_{\frak{Alg}(\Md +E )} $. Then $\Md+E \vdash \Inv_P$. Hence $(\mathbb{Z}/p\mathbb{Z})_0 \not \models E$ for all primes $p$. For the converse, assume $(\mathbb{Z}/p\mathbb{Z})_0 \not \models E$ for all primes $p$. Since $\mathcal{I}_{\frak{Alg}(\Md +E )}$ is minimal, it is a subdirect product of minimal cancellation meadows by Theorem \ref{cancel}. Hence by nontriviality, $\mathcal{I}_{\frak{Alg}(\Md +E )}$ is a subdirect product of  prime expansion fields  $(\mathbb{Z}/p\mathbb{Z})_0$ and $\mathbb{Q}_0$. Since every  $\mathbb{Z}/p\mathbb{Z}_0$ is not a model of $E$, it follows that $\mathcal{I}_{\frak{Alg}(\Md +E )} $ must be isomorphic to $\mathbb{Q}_0$. \hfill $\Box$

It follows that $\Md + \Inv_p$ is the weakest initial algebra specification of $\mathbb{Q}_0$. 
\begin{theorem}
$\frak{Alg}(\Md + \Inv_P)$ is the largest subvariety of $\frak{Alg}(\Md)$ with initial algebra $\mathbb{Q}_0$. 
\end{theorem}
{\bf Proof}: 
By Theorem \ref{initialspec}, $\frak{Alg}(\Md + \Inv_P)$ is a subvariety of $\frak{Alg}(\Md)$ with initial algebra $\mathbb{Q}_0$. Moreover, if $V$ is a subvariety of $\frak{Alg}(\Md)$ with initial algebra $\mathbb{Q}_0$, then every algebra $\mathcal{M}\in V$ must satisfy $\Inv_P$. Thus
$V\subseteq \frak{Alg}(\Md + \Inv_P)$. 
\hfill $\Box$

Finite  expansion fields do not satisfy $\Inv_P$; therefore $\frak{Alg}(\Md + \IL)\not \subseteq \frak{Alg}(\Md + \Inv_P)$. Conversely, not every meadow in $\frak{Alg}(\Md + \Inv_P)$ is a cancellation meadow.
\begin{theorem}
$\frak{Alg}(\Md + \Inv_P)\not \subseteq \frak{Alg}(\Md + \IL)$
\end{theorem}
{\bf Proof}: 
Choose a new constant symbol $a$. For $k\in \mathbb{N}^{+}$ let
\[
E_k = \{a \neq 0\} \cup \{\underline{n}\neq 0 \mid 0 < n < k\} \cup \{\underline{n}\cdot \underline{n}^{-1}=1 \mid 0 < n < k\} \cup  \{\underline{2}\cdot a = 0\} \cup \Md.
\]
Moreover, choose a prime $p \neq 2$ exceeding $k$ and interpret $a$ in $\mathbb{Z}/2p\mathbb{Z}$ by $p$. Since $2p$ is squarefree, $\mathbb{Z}/2p\mathbb{Z}$ is meadow (see \cite{BRS15}). Thus $\mathbb{Z}/2p\mathbb{Z}\models E_k$. 
It follows that $E_k$ is consistent and therefore $E = \bigcup_{k=1}^\infty E_k$ 
is consistent by the compactness
theorem. Let $\mathcal{M}$ be a model for $E$. Then
\begin{enumerate}
\item $\mathcal{M}$ is a meadow, since $\Md \subseteq E$,
\item $\mathcal{M}$ satisfies $\Inv_P$,
\item $\mathcal{M}$ is not a cancellation meadow: $\mathcal{M}\models  \underline{2} \neq 0$ but $\mathcal{M} \models \underline{2}\cdot \underline{2}^{-1}\neq  1$, for otherwise
\[
a = 1 \cdot a =  \underline{2}\cdot \underline{2}^{-1} \cdot a = \underline{2}^{-1}\cdot \underline{2} \cdot a= 2\cdot 0 = 0. \hfill \Box
\]
\end{enumerate}

A frequently asked question of universal algebra is whether or not the identities valid in a class of algebras are \emph{finitely based}, i.e.\ are the consequence of a finite number of identities. Below we  give a negative answer in the case $\Inv_P$.
\begin{theorem}\label{nofin}
$\frak{Alg}(\Md + \Inv_P)$ has no finite basis.
\end{theorem}
{\bf Proof}: Suppose $\frak{Alg}(\Md + \Inv_P)$ has a finite base. Then by the compactness theorem for equational logic, 
there must exist a finite set $R\subseteq P$ such that $\Md+ \{\underline{p}\cdot \underline{p}^{-1} \mid p\in R\}\vdash \Md + \Inv_P$. In order to obtain a contradiction, we choose a prime $p\in P$ larger then any prime in $R$. Then $(\mathbb{Z}/p\mathbb{Z})_0 \in \frak{Alg}(\Md +\{\underline{p}\cdot \underline{p}^{-1} \mid p\in R\} )$ but $(\mathbb{Z}/p\mathbb{Z})_0\not \models \underline{p}\cdot \underline{p}^{-1}=1$. \hfill $\Box$

{\mbox{} It follows  that $\mathbb{Q}_0$ cannot have a finite complete axiomatization which also holds in $\mathbb{C}_0$|the expansion field of the complex numbers. This is proved in Corollary \ref{nocom} below.
\begin{theorem}\label{notc}
Let $s,t$ be  meadow terms. Then
\[
\mathbb{C}_0 \models s=t \ \ \ \text{ iff } \ \ \ \frak{Alg}(\Md + \Inv_P)\models s=t.
\]
\end{theorem}
{\bf Proof}: ($\Leftarrow$): This follows from the fact that $\mathbb{C}_0$ is a cancellation meadow.\\
($\Rightarrow$): Suppose $\frak{Alg}(\Md + \Inv_P)\not \models s=t$. Then there exists a meadow $\mathcal{M}\in \frak{Alg}(\Md  +\Inv_P)$ with $\mathcal{M}\not \models s=t$. Since $\mathcal{M}$ is a subdirect product of cancellation meadows, there exists a cancellation meadow $\mathcal{M}'\in  \frak{Alg}(\Md  + \IL+\Inv_P)$ with $\mathcal{M}'\not \models s=t$. 
Observe that $\mathcal{M}'$ is a cancellation meadow of characteristic 0. Let $\hat{\mathcal{M}'}$ be the algebraic closure of $\mathcal{M}'$. Since $s=t$ is a universal proposition, we have $\hat{\mathcal{M}'}\not \models s=t$. In \cite{BBP15} it is proved that every meadow equation has a first-order representation over the signature $\{0,1, +, \cdot, -\}$ of fields. In particular, there exists a quantifier-free first-order formula $\phi(s,t)$ such that 
\[
\Md + \IL \vdash s=t \leftrightarrow \phi(s,t).
\]
Thus 
\[
\hat{\mathcal{M}'}\models \exists \vec{x}\  \neg \phi(s,t)
\]
where $\exists \vec{x}\  \neg \phi(s,t)$ is the existential closure of $\phi(s,t)$. Since $\mathcal{M}'$ has characteristic 0, $\hat{\mathbb{Q}}_0$ can be embedded in $\hat{\mathcal{M}'}$. Moreover, since algebraically closed fields are model-complete, 
every embedding between algebraically closed fields is elementary (see e.g.\ \cite{R77}). 
We therefore have $\hat{\mathbb{Q}}_0\models \exists \vec{x}\  \neg \phi(s,t)$ and hence $\mathbb{C}_0\models \exists \vec{x}\  \neg \phi(s,t)$. It follows that  $\mathbb{C}_0 \not \models s=t$. \hfill $\Box$
\begin{corollary}\label{nocom}
Let $V\subseteq \frak{Alg}(\Md)$ be a finitely based variety. Then
\[
\mathcal{I}_V \cong \mathbb{Q}_0 \Longrightarrow \mathbb{C}_0 \not \in V.
\]
\end{corollary}
{\bf Proof}: Put $V=\frak{Alg}(\Md + E)$ for some finite set of equations $E$. Suppose $\mathbb{C}_0\models E$. Then $\Md + \Inv_P \vdash E$ by the previous theorem, and hence by equational compactness, $\Md + \{\underline{p}\cdot \underline{p}^{-1}=1\mid p\in R \}\vdash E$ for some finite set $R\subset P$. It follows that $\frak{Alg}(\Md + \Inv_P)$ is finitely based, a contradiction with Theorem \ref{nofin}. \hfill $\Box$}

There exist finitely based subvarieties  $V\subseteq \frak{Alg}(\Md + \Inv_P)$ with initial algebra $\mathbb{Q}_0$. E.g.\  in \cite{BM2011} it is shown that  $\Md + \{(1+x_1^2 + x_2^2)\cdot (1+x_1^2 + x_2^2 )^{-1}=1\}$  is an initial algebra specification of $\mathbb{Q}_0$. There also exist finite initial algebra specifications of $\mathbb{Q}_0$  of the form $\Md+ e$ where $e$ is a single variable equation.
\begin{theorem}\label{single}
Let $p_0$ and $p_1$ be primes and put $f(x)= 2\cdot (x^2 -p_0)(x^2-p_1)(x^2-p_0\cdot p_1)$. Then
\[
\mathbb{Q}_0 \cong \mathcal{I}_{\frak{Alg}(\Md +\{f(x)\cdot f(x)^{-1}=1\})}.
\]
\end{theorem}
{\bf Proof}:
Given a prime $p$ and a natural number $0<n<p$, $n$ is called a {\em quadratic residue} of $p$ if there exists a natural number $x$ such that $x^2\equiv n \mod p$. If the congruence is insoluble, $n$ is said to be a {\em quadratic non-residue}. The quadratic residues and non-residues of an odd prime have a simple \emph{multiplicative property}: the product of two residues or of two non-residues is a residue.

Now consider the polynomial $f(x)= 2 \cdot (x^2-p_0)(x^2-q_0)(x^2-p_0\cdot q_0)$. Inspection shows that $f(x)$ has no rational root. Thus $\mathbb{Q}_0\models f(x)\cdot  f(x)^{-1}=1$. On the other hand, $f(x)\cdot f(x)^{-1}\neq 1$ in every finite prime field $(\mathbb{Z}/p\mathbb{Z})_0$. If $p=2$, this is immediate.
If $p>2$ we apply the multiplicative property of quadratic residues and non-residues. If $(x^2 - p_0)$ or $(x^2 -p_1)$ have a root modulo $p$, then $f(x)$ has a root modulo $p$. If neither $(x^2-p_0)$ nor $(x^2-p_1)$ has  a root modulo $p$ then both $p_0$ and $p_1$ are non-residues of $p$. Hence $p_0\cdot p_1$ is a residue of $p$, i.e.\ $(x^2-p_0\cdot p_1)$ has a root modulo $p$ and thus $f(x)$ has a root modulo $p$. Thus $(\mathbb{Z}/p\mathbb{Z})_0\not \models f(x)\cdot f(x)^{-1}=1$. By Theorem \ref{initialspec} we may conclude that $\Md + f(x)\cdot  f(x)^{-1}=1$ is an initial algebra specification of $\mathbb{Q}_0$. \hfill $\Box$

\section{Expansions of meadows}
In this section we consider expansions of meadows with new function symbols. If $\Sigma$ is a set  of new function symbols  and $E$ is a set of identities, we shall write $\frak{Alg}_\Sigma (\Md +E)$ for the variety  of $\Sigma$-expansion meadows that satisfy $E$; likewise, we shall write  $\frak{Alg}_\Sigma (\Md + \IL +E)$ for the $\Sigma$-expansions of cancellation meadows that satisfy $E$. Thus if  $\Sigma = \emptyset = E$ we have the variety of meadows and the class of cancellation meadows. 
Moreover, if $e$ is an identity we shall write $\frak{Alg}_\Sigma (\Md +E)\models e$ ($\frak{Alg}_\Sigma (\Md +\IL + E)\models e$) if all meadows in $\frak{Alg}_\Sigma (\Md +E)$ ($\frak{Alg}_\Sigma (\Md +\IL + E)$) satisfy $e$, and we shall write $\frak{Alg}_\Sigma (\Md +E)\models E'$ ($\frak{Alg}_\Sigma (\Md +\IL + E)\models E'$) if $\frak{Alg}_\Sigma (\Md +E)\models e$ ($\frak{Alg}_\Sigma (\Md +\IL + E)\models e$) for all $e \in E'$.
Finally, if $\mathcal{A}$ is a $\Sigma$-expansion meadow we shall denote the reduct|the underlying meadow|by $\mathcal{A}\upharpoonright_\Md$. 

Since every algebra is a subdirect product of subdirectly irreducible algebras, every variety $V$ is  determined by its subdirectly irreducible members. 
The following theorem is a semantic version of the \emph{Generic Basis Theorem} proved in \cite{BBP13}.
\begin{theorem}\label{subdir}
Let $\Sigma$ be a set of new function symbols and let $E$ be a set of identities. If every subdirectly irreducible algebra in $\frak{Alg}_\Sigma (\Md +E)$ is a cancellation meadow, then
\[
\frak{Alg}_\Sigma (\Md +E)\models s=t\ \ \ \text{ iff }\ \ \  \frak{Alg}_\Sigma (\Md +\IL + E)\models s=t
\]
for all $\Sigma$-meadow terms $s,t$.
\end{theorem}
{\bf Proof}: Suppose every subdirectly irreducible meadow in  $\frak{Alg}_\Sigma (\Md +E)$ is a cancellation meadow.
Clearly, if $\frak{Alg}_\Sigma (\Md +E)\models s=t$ then $\frak{Alg}_\Sigma (\Md +\IL + E)\models s=t$ for all $s,t$. For the converse,
suppose $\frak{Alg}_\Sigma (\Md +\IL + E)\models s=t$ and let $\mathcal{A}\in \frak{Alg}_\Sigma (\Md +E)$. 
Then  $\mathcal{A}$ is a subdirect product of a family $(\mathcal{A}_i)_{i\in I}$ of subdirectly irreducible meadows in $\frak{Alg}_\Sigma (\Md +E)$. Moreover, since $\mathcal{A}_i\in \frak{Alg}_\Sigma (\Md +\IL + E)$ for every $i$, it follows that 
$\mathcal{A}_i$ satisfies $s=t$ for every $i$. Hence $\prod_{i\in I} \mathcal{A}_i$ satisfies $s=t$. Thus $\mathcal{A}$ satisfies $s=t$, since $\mathcal{A}$ is a subalgebra of $\prod_{i\in I} \mathcal{A}_i$. \hfill $\Box$

A typical expansion variety with non-cancellation subdirectly irreducible meadows is the following. Consider the expansion with a single binary function symbol $\mathbf{eq}$ and the set $E$ of equations consisting of $\mathbf{eq}(x,x)=1$ and $\mathbf{eq} (\underline{i},\underline{j})=0$ for $0\leq i\neq j\leq 5$ which define equality in $\mathbb{Z}/6\mathbb{Z}$. One can interpret $\mathbf {eq}$ in the meadow $\mathcal{M}= (\mathbb{Z}/2\mathbb{Z})_0 \times (\mathbb{Z}/3\mathbb{Z})_0$ by
\[
\mathbf{eq}^\mathcal{M}(\langle i,i'\rangle, \langle j,j'\rangle) =
\begin{cases}
\langle 1,1\rangle & \text{ if  $i=j$ and $i'=j'$ },\\
\langle 0,0\rangle  & \text{ otherwise.}
\end{cases}
\]

$\mathcal{M}\upharpoonright_\Md$ is not subdirectly irreducible: $\Theta_1=ker(\pi_1)$ and $\Theta_2=ker(\pi_2)$ are both congruences but $\Theta_1\cap \Theta_2= \Delta$. $\mathcal{M}$, however, is subdirectly irreducible: its only congruences are the diagonal and the all relation. Moreover, $\mathcal{M}$ is not  a cancellation meadow: since
\[
\langle1,0\rangle \cdot \langle1,0\rangle \cdot \langle1,0\rangle = \langle1,0\rangle 
\]
it follows that $\langle1,0\rangle^{-1}= \langle1,0\rangle$. But 
\[
\langle1,0\rangle \cdot \langle1,0\rangle ^{-1}= \langle1,0\rangle \cdot \langle1,0\rangle= \langle1,0\rangle \neq \langle1,1\rangle .
\]
Also, e.g.\ $\frak{Alg}_\Sigma (\Md +\IL + E)\models \mathbf{eq}(x\cdot x^{-1},1)=x\cdot x^{-1}$; but $\frak{Alg}_\Sigma (\Md + E)\not \models \mathbf{eq}(x\cdot x^{-1},1)=x\cdot x^{-1}$ since 
\[
\mathbf{eq}^\mathcal{M}(\langle 1,0\rangle \cdot \langle 1,0\rangle^{-1}, \langle 1,1\rangle) =\mathbf{eq}^\mathcal{M}(\langle 1,0\rangle, \langle 1,1\rangle) =\langle 0,0\rangle \neq  \langle 1,0\rangle=\langle 1,0\rangle \cdot \langle 1,0\rangle^{-1}.
\]
If $\Sigma=\emptyset$, we have the following corollary.
\begin{corollary}\label{mdvar}
Let $E$ be a set of meadow equations. Then 
\begin{enumerate}
 \item 
$
\frak{Alg}(\Md +E)\models s=t \text{ iff } \frak{Alg} (\Md +\IL +E)\models s=t $
for all meadow terms $s,t$, and
\item if $E$ is finite then $\frak{Alg}(\Md +E )$  can be axiomatized by at most eleven equations.  
\end{enumerate}
\end{corollary}
{\bf Proof}: 
\begin{enumerate}
\item This follows from Theorem \ref{subdir} and the fact  that every subdirectly irreducible meadow is a cancelation meadow. 
\item We first prove that 
\[
\frak{Alg}(\Md + \{r=0, t=0\})= \frak{Alg}(\Md + (1-t\cdot t^{-1})(1-r \cdot r^{-1})=1). \ \ \ \ (\ddag)
\]
Clearly, $\frak{Alg}(\Md + \{r=0, t=0\})\subseteq \frak{Alg}(\Md + (1-t\cdot t^{-1})(1-r \cdot r^{-1})=1)$. In order to prove the converse inclusion, we have to show 
\[
\frak{Alg}(\Md + (1-t\cdot t^{-1})(1-r \cdot r^{-1})=1)\models  \{r=0, t=0\}.
\]
By  1. it suffices to prove 
\[
\frak{Alg}(\Md + \IL+ (1-t\cdot t^{-1})(1-r \cdot r^{-1})=1)\models  \{r=0, t=0\}.
\]
This follows immediately, because if $\mathcal{M}\in \frak{Alg}(\Md + \IL+ (1-t\cdot t^{-1})(1-r \cdot r^{-1})=1)$ and e.g.\ $\mathcal{M}\not\models t=0$ then there exist $t^*, r^*\in M$ with $t^*\neq 0$ and $(1-t^*\cdot t^{*-1})(1-r^* \cdot r^{*-1})=1$. It follows that $\mathcal{M}\models 0=1$, since 
$ t^*\cdot t^{*-1}=1$.

Now let $E$ be a finite set of equations. We may assume that every equation is of the form $s=0$. By $(\ddag)$ we can code two equations of the form $r=0$ and $t= 0$ in the  equation $(1-t\cdot t^{-1})(1-r \cdot r^{-1}) - 1 =0$. By the same argument  we can proceed reducing smaller sets of equations until only a single equation $e$ is left.   \hfill $\Box$
\end{enumerate}

A particular expansion with the sign function $\mathbf s$ was studied in \cite{BBP15} where a finite axiomatization of formally real meadows was given. We define the sign function in an equational manner by the set $\SA$ of axioms given in Table~\ref{t:a}.
\begin{table}
\centering
\hrule
\begin{align}
\sg(1_x)&=1_x&\tag{S1}\label{S1}\\
\sg(0_x)&=0_x&\tag{S2}\label{S2}\\
\sg(-1)&=-1&\tag{S3}\label{S3}\\
\sg(x^{-1})&=\sg(x)&\tag{S4}\label{S4}\\
\sg(x\cdot y)&=\sg(x)\cdot \sg(y)&\tag{S5}\label{S5}\\
0_{\sg(x)-\sg(y)}\cdot (\sg(x+y)-\sg(x))&=0&\tag{S6}\label{S6}
\end{align}
\hrule
\caption{The set $\SA$ of axioms for the sign function}
\label{t:a}
\end{table}
Here we write $1_x$ for $x\cdot x^{-1}$ and $0_x$ for $1-1_x$.  Clearly, $\frak{Alg}_{\{\sg\}}(\Md + \SA)$ is a variety. But also the class of all reducts of signed meadows is a variety. This means that the sign function has an equational specification in the signature of the meadows.
\begin{theorem} \label{sign}
$\{ \mathcal{M}\upharpoonright_\Md \mid \mathcal{M}\in \frak{Alg}_{\{\sg\}}(\Md + \SA)\}$ is a variety.
\end{theorem}
{\bf Proof}: Let $K= \{ \mathcal{M}\upharpoonright_\Md \mid \mathcal{M}\in \frak{Alg}_{\{\sg\}}(\Md + \SA)\}$. In \cite{BBP13} (Proposition 3.9), it is proved that $K\subseteq \frak{Alg}(\Md + \EFR)$ where $\EFR= \{0_{x_0^2 + \cdots +x_n^2}\cdot x_0=0\mid n\in \mathbb{N}\}$ is an infinite axiomatization of formal realness. In order to prove that $K$ is a variety, it therefore suffices to prove that  $\frak{Alg}(\Md + \EFR)\subseteq K$.

We have to show that every $\mathcal{M}\in \frak{Alg}(\Md + \EFR)$ can be expanded to a signed meadow. Thus pick $\mathcal{M}\in \frak{Alg}(\Md + \EFR)$. Then $\mathcal{M}$ is a subdirect product of a family $(\mathcal{M}_i)_{i\in I}$ of formally real cancellation meadows. Every $\mathcal{M}_i$  can thus be ordered by an ordering $<$ satisfying the axioms (OF1)--(OF4) given below: 
\[
\begin{array}{rlr}
x\neq 0 &\rightarrow( x<0 \ \vee\ 0<x )& \text{(OF1)}\\
x < y& \rightarrow\neg ( y < x \ \vee\ x=y)& \text{(OF2)}\\
x<y & \rightarrow  x+z < y+z &\text{(OF3)}\\
x<y \ \wedge\ 0<z & \rightarrow  x\cdot z <y \cdot z&\text{(OF4)}
\end{array}
\]
W.l.o.g.\ we may assume that the ordering is such that $-1 < 0  <1$.  We define the sign function $\sg_i$ on the carrier $M_i$ of $\mathcal{M}_i$ by
\[\sg_i(a)=\begin{cases}
-1&\text{if }a<0,\\
0&\text{if }a=0,\\
1&\text{if }a>0.
\end{cases}
\]
Then $\sg_i$ satisfies the axioms of $\SA$: (S1) and (S2) hold since  $1_a, 0_a \in \{0,1\}$ for every $a\in M_i$.
(S3) is immediate. To prove (S4), observe that  $a=0$ or $a<0$ or $a>0$ by (OF1). Moroever, $a=0$ iff $a^{-1}=0$ and, by (OF4), $a<0$ iff $a^{-1}<0$. Hence also (S4) holds. (S5) follows from a case distinction.  Finally, (S6) is an equational representation of the conditional axiom $\sg(x)=\sg(y) \longrightarrow \sg(x+y)=\sg(x)$ which clearly holds too.

Since $\mathcal{M}$ is a subalgebra of $\prod_{i\in I} \mathcal{M}_i$, we can now define the sign function on $\mathcal{M}$ in the obvious way by $\sg(\langle a_i\rangle_{i \in I})= \langle \sg_i(a_i)\rangle_{i \in I}$. Then $\mathcal{M}$ is a signed meadow. \hfill $\Box$

The situation is different if we expand signed meadows with a \emph{signed square root} $\surd$ defined by the equations given in Table \ref{surd}.
\begin{table}\label{surd}
\centering
\hrule
\begin{align}
\surd (x^{-1})&=(\surd (x))^{-1}\\
\surd (x\cdot y)&=\surd (x)\cdot \surd (y)\\
\surd (x\cdot x \cdot \sg (x))&=x\\
\sg ( \surd (x) - \surd (y)) &= \sg (x - y)
\end{align}
\hrule
\caption{The set $SR$ of axioms for the signed square root}
\end{table}
Here we stipulate $\surd (x) = - \surd (-x)$ for $x < 0$. The class of reducts of signed meadows with square roots $K=\{ \mathcal{M}\upharpoonright_\Md \mid \mathcal{M}\in \frak{Alg}_{\{\sg, \surd \}}(\Md + \SA + \SR)\}$ is not a variety: clearly $\mathbb{R}_0\in K$ but its subalgebra $\mathbb{Q}_0$ cannot be expanded with square roots. Hence $K$ is not closed under subalgebras. It follows that unlike the sign function, the square roots have no equational specification in the signature of the meadows.

\section{Presentations and specializations}
A finite algebra may be presented by a table showing all its elements and the values of its operations. In the infinite case we employ an initial algebra.
\begin{definition}
A \emph{presentation} of $\mathcal{M} \in \frak{Alg}(\Md)$ is a pair $(C,E)$ of fresh constants and identities such that the following holds:
\begin{enumerate}
\item $E$ is a set of meadow equations over $C$, and
\item $\mathcal{M} \cong \mathcal{I}_{\frak{Alg}_{C}(\Md +E)}\upharpoonright_\Md$.
\end{enumerate}
A presentation $(C,E)$ of $\mathcal{M}$ is \emph{finite} if both $C$ and $E$ are finite.
\end{definition}
The meadow of Gaussian rationals|denoted $\mathbb{Q}_0(i)$| is obtained by adjoining the imaginary number $i$ to the meadow of rationals. In \cite{BB16} it is shown that $(\{i\}, \{f(x)\cdot f(x)^{-1}=1, i^2 +1=0\})$ where $f(x)=(x^2-2)(x^2-3)(x^2-6)$ is a finite presentation of $\mathbb{Q}_0(i)$. The following theorem is a generalization of this fact.
\begin{theorem}\label{finpres}
Let $\mathbb{Q}(c_1, \ldots , c_n)$ be an algebraic extension of $\mathbb{Q}$. Then the meadow $\mathbb{Q}_0(c_1, \ldots , c_n)$ has a finite presentation.
\end{theorem}
{\bf Proof}: Recall the {\em Primitive element theorem} which says that if $E$ is a finite degree separable extension over field $F$ then $E=F(\alpha)$ for some $\alpha \in E$. This theorem applies to algebraic number fields, i.e.\ finite extensions of the rational numbers 
$\mathbb{Q}$, since $\mathbb{Q}$ has characteristics $0$ and therefore every  extension over $\mathbb{Q}$ is separable.

By the {\em Primitive element theorem} we may assume that $\mathbb{Q}(c_1, \ldots , c_n)= \mathbb{Q}(c)$ for some well-chosen $c$. Let $g(x)= q_0 + q_1\cdot x + \cdots + q_n \cdot x^n$ be the minimal polynomial of $c$ over $\mathbb{Q}$. Since $g$ is minimal it follows that it is monic, i.e.\ $q_n=1$. Now choose primes $p$ and $q$ such that the polynomial $f(x)=2(x^2 - p)(x^2-q)(x^2-p\cdot q)$
has no root in $\mathbb{Q}(c)$. That  $p$ and $q$ exist can be seen as follows. Since $\mathbb{Q}(c)$ is a finite extension over $\mathbb{Q}$, it contains only finitely many subfields. And if $p$ and $q$ are different primes, then $\mathbb{Q}(\sqrt{p})$, $\mathbb{Q}(\sqrt{q})$ and $\mathbb{Q}(\sqrt{p\cdot q})$ are non-isomorphic subfields of $\mathbb{Q}(c)$. So $p$ and $q$ exist. Now consider the pair $\mathcal{P}=(\{c\}, \{f(x)\cdot f(x)^{-1}=1, g(c)=0\})$. We show that $\mathcal{P}$ is a presentation of $\mathbb{Q}(c)$. Clearly 
$\mathbb{Q}_0(c)\models \Md + \{f(x)\cdot f(x)^{-1}=1, g(c)=0\}$. Moreover by Theorem \ref{single}, $\mathbb{Q}_0 \cong \mathcal{I}_{\frak{Alg}(\Md +\{f(x)\cdot f(x)^{-1}=1\})}$.
It therefore suffices to show that every closed meadow term over $c$ is provable equal to a term of the form $a_0 + a_1\cdot c + \cdots + a_{n-1}\cdot c^{n-1}$ where every $a_i$ is of the form $\underline{n}\cdot \underline{m}^{-1}$. We employ structural induction.

Clearly, this holds for $0$, $1$ and $c$. For $-t$, $t +t'$ and $t\cdot t'$ it holds as well, since the equation $g(c)=0$ can be used to eliminate all powers of $c$ higher than $n-1$. Now consider $t^{-1}$. By the induction hypothesis there exists an appropriate polynomial $p(x)=a_0 + a_1\cdot x + \cdots  + a_{n-1}\cdot x^{n-1}$ such that  $t^{-1}= p(c)^{-1}$. If $a_i=0$ for all $i$, we are done. Thus assume that not all $a_i$ are $0$. Then  $gcd(p(x), g(x))=1$, since $g$ is minimal.  It follows that  the Eucledian division procedure of polynomials produces appropriate polynomials
$h(x)$ and $h'(x)$ such that $h(x)\cdot g(x) + h'(x)\cdot p(x) =1$. Hence $h'(c)\cdot p(c)=1$ since $g(c)=0$. So $t^{-1}= h'(c)$. \hfill $\Box$

Classes of algebras which are not varieties are often introduced as specialized subclasses of varieties. In the remainder of this section we consider specializations of meadow varieties by constants.
\begin{definition}
A \emph{specialization by constants} of  a variety $ \frak{Alg}(\Md +E )\subseteq \frak{Alg}(\Md)$ is a pair $(C,E')$ of fresh constants and identities such that the following holds:
\begin{enumerate}
\item $E'$ is a set of meadow equations over $C$, and
\item $\{\mathcal{M} \upharpoonright_\Md\mid \mathcal{M} \in \frak{Alg}_{C}(\Md +E +E' )\}\subseteq \frak{Alg}(\Md +E )$.
\end{enumerate}
A specialization  of $\frak{Alg}(\Md +E )$ is \emph{finite} if both $C$ and $E'$ are finite, is \emph{stricly smaller} if $\{\mathcal{M} \upharpoonright_\Md\mid \mathcal{M} \in \frak{Alg}_{C}(\Md +E +E' )\}\neq \frak{Alg}(\Md +E )$ and is nontrivial if $\frak{Alg}_{C}(\Md +E +E')$ is nontrivial. 
\end{definition}
E.g. $\frak{Alg}_{\{\sqrt{2}\}}(\Md +  \{\sqrt{2}\cdot \sqrt{2}=\underline{2}\})$ is a finite, strictly smaller, nontrivial
specialization of $\frak{Alg}(\Md)$ since $\mathbb{Q}_0$ has no expansion in $\frak{Alg}_{\{\sqrt{2}\}}(\Md +\{\sqrt{2}\cdot \sqrt{2}=\underline{2}\})$.
\begin{proposition}
Every nontrivial variety $V \subseteq  \frak{Alg}(\Md)$ has a nontrivial, strictly smaller specialization by constants.
\end{proposition}
{\bf Proof}:
Assume $V= \frak{Alg}(\Md +E)$ is nontrivial, say $V$ contains a meadow $\mathcal{M}$ of cardinality $\kappa >1$. Let $C$ be a set containing  $\kappa^+$ fresh constants and add to the equations  in $E$ the set $E'=
\{(c- c')\cdot (c- c')^{-1} =1\mid  c\neq c' \in C\}$. Then $(C, E')$ is a nontrivial specialization of $V$. Moreover, $(C, E')$ is strictly smaller since $\mathcal{M}$ has no expansion  in
$\frak{Alg}_C(\Md +E +E')$.  \hfill $\Box$

Clearly, if $V$ contains a finite, nontrivial meadow, then it also has a finite, nontrivial specialization which is strictly smaller. On the other hand, if $V\models \Inv_P$ then there always exists a specialization which has an initial element that cannot be finitely presented.
\begin{theorem}
Let $V\subseteq  \frak{Alg}(\Md + \Inv_P)$ be a nontrivial variety. Then $V$ has a specialization by constants  $ \frak{Alg}_{\{c\}}(\Md +E)$ such that 
\begin{enumerate}
\item $\mathcal{I}_{ \frak{Alg}_{\{c\}}(\Md +E)}$ is generated by $c$, and 
\item $\mathcal{I}_{ \frak{Alg}_{\{c\}}(\Md +E)}\upharpoonright_\Md$ has no finite presentation.
\end{enumerate}
\end{theorem}
{\bf Proof}: Suppose $V=\frak{Alg}(\Md +E')$. For $U\subseteq \mathbb{N}$ define $\Gamma_U=\{(c-\underline{n})\cdot (c-\underline{n})^{-1}=1 \mid n\in U\}$. Put $E_U=\Md + E' + \Gamma_U$ and $V_U=\frak{Alg}_{\{c\}}(E_U)$. We first prove that 
\[
(\dag)\ \ \ \ \mathcal{I}_{V_U} \models (c-\underline{n})\cdot (c-\underline{n})^{-1}=1\ \ \ \  \text{ iff }\ \ \ \  n\in U. 
\]
Right-to-left is immediate. For the converse assume $n\not\in U$. Choose a nontrivial meadow $\mathcal{M}\in V$ and interpret $c$ in $\mathcal{M}$ by $\underline{n}$. Then $\mathcal{M}\in V_U$, since $\mathcal{M}\models \Inv_P$, and $\mathcal{M}\not\models (c-\underline{n})\cdot (c-\underline{n})^{-1}=1$. Hence $\mathcal{I}_{V_U} \not\models (c-\underline{n})\cdot (c-\underline{n})^{-1}=1$.

Now assume that for given different $V,W \subseteq \mathbb{N}$ we have $\mathcal{I}_{V_U}\upharpoonright_\Md \cong \mathcal{I}_{V_W}\upharpoonright_\Md$. Then there exists an isomorphism $\phi: \mathcal{I}_{V_U}\upharpoonright_\Md \rightarrow \mathcal{I}_{V_W}\upharpoonright_\Md$. Pick a closed meadow term $t$ over $c$ with $\phi([c])=[t]$. Then
\[
\begin{array}{rcll}
n\in U & \longleftrightarrow & \mathcal{I}_{V_U} \models (c-\underline{n})\cdot (c-\underline{n})^{-1}=1 & \text{by $\dag$,}\\
&  \longleftrightarrow & \mathcal{I}_{V_W} \models (t-\underline{n})\cdot (t-\underline{n})^{-1}=1 & \text{since $\phi$ is a homomorphism.}
\end{array}
\]
Hence $U$ is recursively enumerable in $W$; likewise $W$ is recursively enumerable  in $U$. Now let $\equiv$ be the equivalence relation on $\mathcal{P}(\mathbb{N})$ defined by 
\[
U \equiv W \Longleftrightarrow \text{ $U$ is recursively enumerable in $W$ and $W$ is recursively enumerable  in $U$.}
\]
Then every equivalence class is countable; consequently there are $2^{\aleph_0}$ equivalence classes.
However, the number of equivalence classes with finite presentations is countable. Hence there must exist a specialization by constants $\frak{Alg}_{\{c\}}(\Md +E)$ such that $\mathcal{I}_{ \frak{Alg}_{\{c\}}(\Md +E)}\upharpoonright_\Md$ has no finite presentation. \hfill $\Box$

\section{Open questions}
We suspect that the question whether a finite axiomatization of $\mathbb{Q}_0$ exists can be solved if the subvarieties of the meadows were well understood and fully classified. We end this paper with some open questions which we believe can approach the problem from different angles.

Finitely based initial algebra specifications of $\mathbb{Q}_0$ exist in various forms:
e.g.\  in \cite{BM2011} it is shown that  $\frak{Alg}(\Md + (1+x_1^2 + x_2^2)\cdot (1+x_1^2 + x_2^2 )^{-1}=1)$ is a subvariety of meadows with initial element $\mathbb{Q}_0$; in Theorem \ref{single} we proved that a single one-element identity can be added to $\Md$ in order to yield an initial algebra specification.
The question then arises whether there exist  maximal or minimal finitely based subvarieties of $\frak{Alg}(\Md + \Inv_P)$.   If there exists a minimal finitely  based subvariety of $\frak{Alg}(\Md + \Inv_P)$, then a finite complete axiomatization of the 
equational theory of the rational numbers exists. In particular, by Corollary \ref{mdvar} a complete axiomatization of the form $\Md + e$ can be given. It is unclear whether $e$ can be a one-variable equation.

Another way to tackle the problem is to study finite presentations of $\mathbb{Q}_0$ and finite specializations of $\frak{Alg}(\Md + \Inv_P)$. In Theorem \ref{finpres} we showed that every finite algebraic extension of $\mathbb{Q}_0$ has a finite presentation. Does this fact hold for arbitrary  extension fields of $\mathbb{Q}_0$, i.e.\ extension fields that have  at least one element that is transcendental over $\mathbb{Q}_0$?

{\bf Acknowledgement:} We acknowledge a useful remark by Gerard van der Geer (Korteweg--de Vries Institute, University of Amsterdam) concerning the proof of Theorem 12.

\end{document}